\documentclass[12pt,a4paper]{article}
\usepackage[T1]{fontenc}
\usepackage[scaled=0.91]{helvet}
\usepackage{courier}
\usepackage{amsmath,amssymb,amsfonts}        % moduli per le estensioni ed i fonts AMS
\usepackage{mathptmx}
\usepackage[dvips]{graphicx}
\usepackage{psfrag}
\usepackage{color}
\usepackage[pdfauthor={Riccardo Fazio},pdftitle={Error estimator for FD}, bookmarks,colorlinks]{hyperref}

\def\ds{\displaystyle}
\def\RE{Richardson's extrapolation}
\def\FBF{free boundary formulation}
\def\www{.9\textwidth}
\def\RR{\hbox{I\kern-.2em\hbox{R}}}
\def\NN{\hbox{I\kern-.2em\hbox{N}}}
\def\ds{\displaystyle}

\def\E{\mbox{E}}

\newcommand{\eqnsection}{
   \renewcommand{\theequation}{{\thesection.\arabic{equation}}}
   \makeatletter
   \csname @addtoreset\endcsname{equation}{section}
   \makeatother}
\eqnsection

\title{A Posteriori Error Estimator for a\\
Non-Standard Finite Difference Scheme\\
Applied to BVPs on Infinite Intervals}
\author{Riccardo Fazio\thanks{Corresponding author home-page: http://mat521.unime.it/$\sim$fazio} \ \ and Alessandra Jannelli \\
Department of Mathematics and Computer Science \\
University of Messina \\
Viale F. Stagno D'Alcontres 31, 98166 Messina, Italy \\
rfazio@unime.it\ \ \ \ \ ajannelli@unime.it}
\date{\today}
\begin{document}
\maketitle
\begin{abstract}
In this paper, we present a study of an a posteriori estimator for the discretization error of a non-standard finite difference scheme applied to boundary value problems defined on an infinite interval.
In particular, we show how Richardson's extrapolation can be used to improve the numerical solution involving the order of accuracy and numerical solutions from two nested quasi-uniform grids. 
A benchmark problem is examined for which the exact solution is known and we get the following result: if the round-off error is negligible and the grids are sufficiently fine then the Richardson's error estimate gives an upper bound of the global error. 
\end{abstract}
\smallskip

\smallskip


%\newpage

\section{Introduction}
The main aim of this paper is to show how Richardson's extrapolation can be used to define an error estimator for a non-standard finite difference scheme applied to boundary value problems (BVPs) defined on the infinite interval.
%\section{BVPs on the half line}
Without loss of generality, we consider the class of BVPs 
\begin{eqnarray}
&& {\ds \frac{d{\bf u}}{dx}} = {\bf f} \left(x, {\bf u}\right)
\ , \quad x \in [0, \infty) \ , \nonumber \\[-1.5ex]
\label{p} \\[-1.5ex]
&& {\bf g} \left( {\bf u}(0), {\bf u} (\infty) \right) = {\bf 0}
\ ,  \nonumber
\end{eqnarray}
where $ {\bf u}(x) $ is a $ d-$dimensional vector with $ {}^{\ell} u
(x) $ for $ \ell =1, \dots , d $ as components, $ {\bf f}:[0,
\infty) \times \RR^d \rightarrow~\RR^d $, and $ {\bf g}:
\RR^d \times \RR^d \rightarrow \RR^d $.
Here, and in the following, we use Lambert's notation for the vector components \cite[pp. 1-5]{Lambert}.
Existence and uniqueness results, as well as results concerning the solution asymptotic behaviour, for classes of problems belonging to (\ref{p}) have been reported in the literature, see for instance Granas et al. \cite{Granas:BVP:1986}, Countryman and Kannan \cite{Countryman:1994:CNB}, and Agarwal et al. \cite{Agarwal:2002:EAB,Agarwal:2002:IIP,Agarwal:2005:EAB}.

Numerical methods for problems belonging to (\ref{p}) can be classified according to the numerical treatment of the boundary conditions imposed at infinity.
The oldest and simplest treatment is to replace the infinity with a suitable finite value, the so-called truncated boundary.
However, being the simplest approach this has revealed within the decades some drawbacks that suggest not to apply it especially if we have to face a given problem without any clue on its solution behaviour.
Several other treatments have been proposed in the literature to overcome the shortcomings  of the truncated boundary approach.
In this research area they are worth of consideration: the formulation of so-called asymptotic boundary conditions by de Hoog and Weiss \cite{deHoog:1980:ATB}, Lentini and Keller \cite{Lentini:BVP:1980} and Markowich
\cite{Markowich:TAS:1982,Markowich:ABV:1983}; the reformulation of the given problem in a bounded domain as studied first by de Hoog and Weiss and developed more recently by Kitzhofer et al. \cite{Kitzhofer:2007:ENS}; the \FBF \ proposed by Fazio \cite{Fazio:1992:BPF} where the unknown free boundary can be identified with a truncated boundary;  the treatment of the original domain via pseudo-spectral collocation methods, see the book by Boyd \cite{Boyd:2001:CFS} or the review by  Shen and Wang \cite{Shen:SRA:2009} for more details on this topic; and, finally, a non-standard finite difference scheme on a quasi-uniform grid defined on the original domain by Fazio and Jannelli \cite{Fazio:2014:FDS}.

When solving a mathematical problem by numerical methods one of the main concerns is related to the evaluation of the global error.
For instance, Skeel \cite{Skeel:1986:TWE} reported on thirteen strategies to approximate the numerical error.
Here we are interested to show how within Richardson's extrapolation theory we can derive an error estimate.
For any component $U$ of the numerical solution, the global error $e$ can be defined by
\begin{equation}\label{eq:GE}
e = u - U \ ,
\end{equation}
where $u$ is the exact analytical solution component.
Usually, we have several different sources of errors: discretization, round-off, iteration and programming errors.
Discretization errors are due to the replacement of a continuous problem with a discrete one and the related error decreases by reducing the discretization parameters, enlarging the value of $N$, the number of grid points in our case.
Round-off errors are due to the utilization of floating-point arithmetic to implement the algorithms available to solve the discrete problem.
This kind of error usually decreases by using higher precision arithmetic, double or, when available, quadruple precision.
Iteration errors are due to stopping an iteration algorithm that is converging but only as the number of iterations goes to infinity.
Of course, we can reduce this kind of error by requiring more restrictive termination criteria for our iterations, the iterations of Newton's method in the present case.    
Programming errors are beyond the scope of this work, but they can be eliminated or at least reduced by adopting the so-called structured programming.  
When the numerical error is caused prevalently by the discretization error and in the case of smooth enough solutions the discretization error can be decomposed into a sum of powers of the inverse of $N$ 
\begin{equation}\label{eq:asymE}
u = U_{N} + C_0 \left(\frac{1}{N}\right)^{p_0}+ C_1 \left(\frac{1}{N}\right)^{p_1}+ C_2 \left(\frac{1}{N}\right)^{p_2}+ \cdots \ ,
\end{equation}
where $C_0$, $C_1$, $C_2$, $\dots$ are coefficients that depend on $u$ and its derivatives, but are independent on $N$, and $p_0$, $p_1$, $p_2$, $\dots$ are the true orders of the error.
The value of each $p_k$, for $k=0$, $1$, $2$, $\cdots$, is usually a positive integer with $p_0 < p_1 < p_2 < \cdots$ and all together constitute an arithmetic progression of ratio $p_1-p_0$, see Joyce \cite{Joyce:1971:SEP}.
The value of $p_0$ is called the asymptotic order or the order of accuracy of the method or of the numerical solution $U$. 
 %\textcolor{red}{cite the bibliography}

\section{Numerical scheme}\label{S:scheme}
In order to solve a problem in the class (\ref{p}) on the original domain we discuss first quasi-uniform grids maps from a reference finite domain and introduce on the original domain a non-standard finite difference scheme that allows us to impose the given boundary conditions exactly.
 
\subsection{Quasi-uniform grids}\label{SS:quniform}
Let us consider the smooth strict monotone quasi-uniform maps $x = x(\xi)$, the so-called grid generating functions,
see Boyd \cite[pp. 325-326]{Boyd:2001:CFS} or Canuto et al. \cite[p. 96]{Canuto:2006:SMF},
\begin{equation}\label{eq:qu1}
x = -c \cdot \ln (1-\xi) \ ,
\end{equation}
and
\begin{equation}\label{eq:qu2}
x = c \frac{\xi}{1-\xi} \ ,
\end{equation}
where $ \xi \in \left[0, 1\right] $, $ x \in \left[0, \infty\right] $, and $ c > 0 $ is a control parameter.
So that, a family of uniform grids $\xi_n = n/N$ defined on interval $[0, 1]$ generates one parameter family of quasi-uniform grids $x_n = x (\xi_n)$ on the interval $[0, \infty]$.
The two maps (\ref{eq:qu1}) and (\ref{eq:qu2}) are referred as logarithmic and algebraic map, respectively. 
As far as the authors knowledge is concerned, van de Vooren and Dijkstra \cite{vandeVooren:1970:NSS} were the first to use this kind of maps. 
We notice that more than half of the intervals are in the domain with length approximately equal to $c$ and 
$x_{N-1} = c \ln N$ for (\ref{eq:qu1}),
while $ x_{N-1} \approx c N $ for (\ref{eq:qu2}).
For both maps, the equivalent mesh in $x$ is nonuniform with the
most rapid variation occurring with $c \ll x$.
The logarithmic map (\ref{eq:qu1}) gives slightly better resolution near $x = 0$ than the
algebraic map (\ref{eq:qu2}), while the algebraic map gives much better resolution than the
logarithmic map as $x \rightarrow \infty$. 
In fact, it is easily verified that
\[
-c \cdot \ln (1-\xi) < c \frac{\xi}{1-\xi} \ ,
\]
for all $\xi$, but $\xi = 0$, see figure \ref{fig:m1N20} below.

The problem under consideration can be discretized by introducing a uniform grid $ \xi_n $ of $N+1$ nodes in $ \left[0, 1\right] $ with $\xi_0 = 0$ and $ \xi_{n+1} = \xi_n + h $ with $ h = 1/N $, so that $ x_n $ is a quasi-uniform grid in $ \left[0, \infty\right] $. 
The last interval in (\ref{eq:qu1}) and (\ref{eq:qu2}), 
namely $ \left[x_{N-1}, x_N\right] $, is infinite but the point $ x_{N-1/2} $ is finite, because the non integer nodes are defined by 
\[
x_{n+\alpha} = x\left(\xi=\frac{n+\alpha}{N}\right) \ ,
\]
with $ n \in \{0, 1, \dots, N-1\} $ and $ 0 < \alpha < 1 $.
%In this way we have defined also the ghost cell $ x_{-1} $.  
These maps allow us to describe the infinite domain by a finite number of intervals.
The last node of such grid is placed at infinity so right boundary conditions
are taken into account correctly.
Figure \ref{fig:m1N20} shows the two quasi-uniform grids $x=x_n$, $n = 0, 1, \dots , N$ defined by (\ref{eq:qu1}) and by (\ref{eq:qu2}) with $c=10$ and $N$ equal to, from top to bottom, 10, 20 and 40, respectively.
\begin{figure}[!hbt]
\centering
\psfrag{x}[][]{$x$} 
\psfrag{y}[][]{$$} 
\psfrag{y1}[][]{$$} 
\includegraphics[width=.8\textwidth,height=6cm]{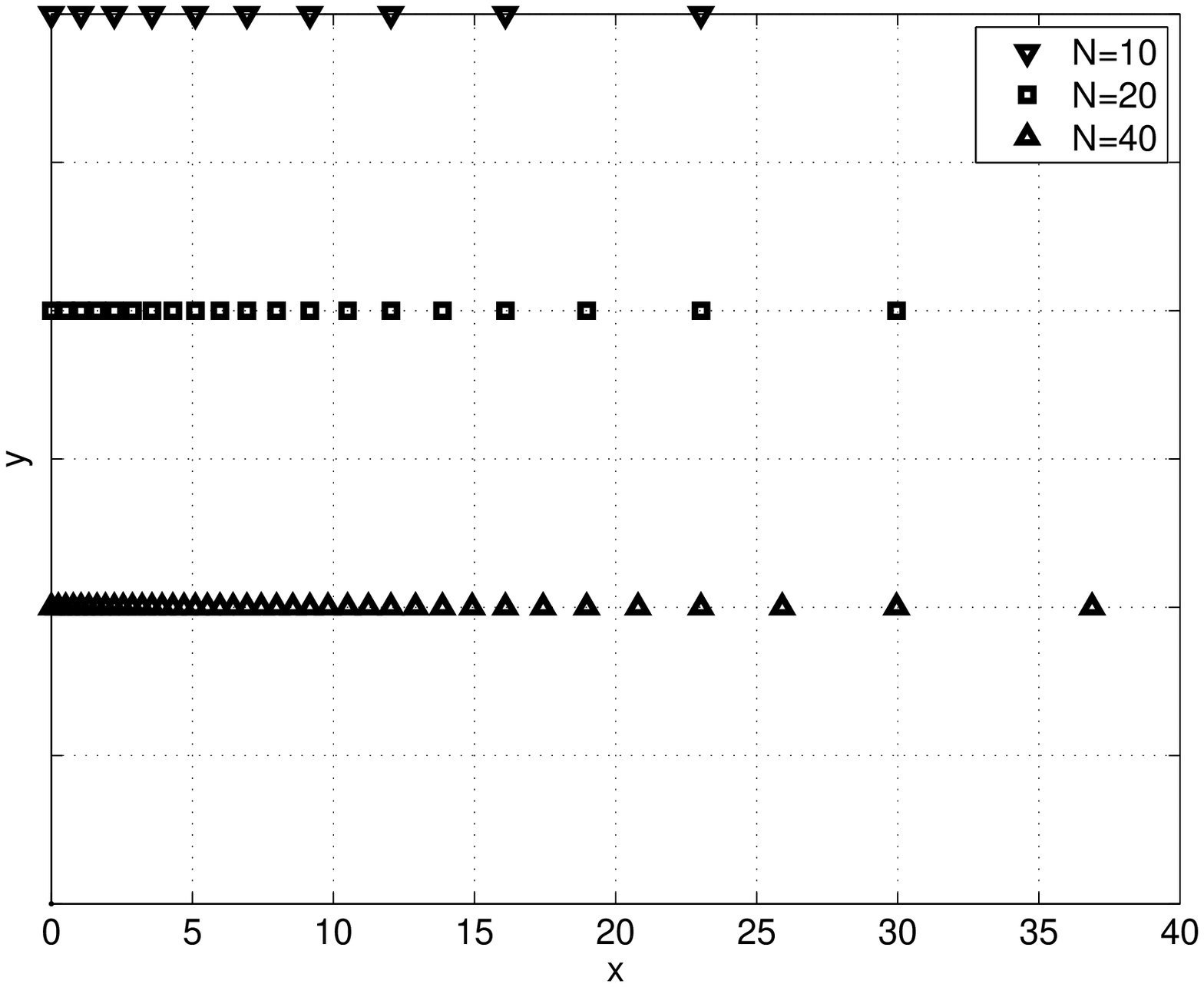} 
\put(5,157){$x_N = \infty$} 
\put(5,108){$x_N = \infty$} 
\put(5,60){$x_N = \infty$} \\
\includegraphics[width=.8\textwidth,height=6cm]{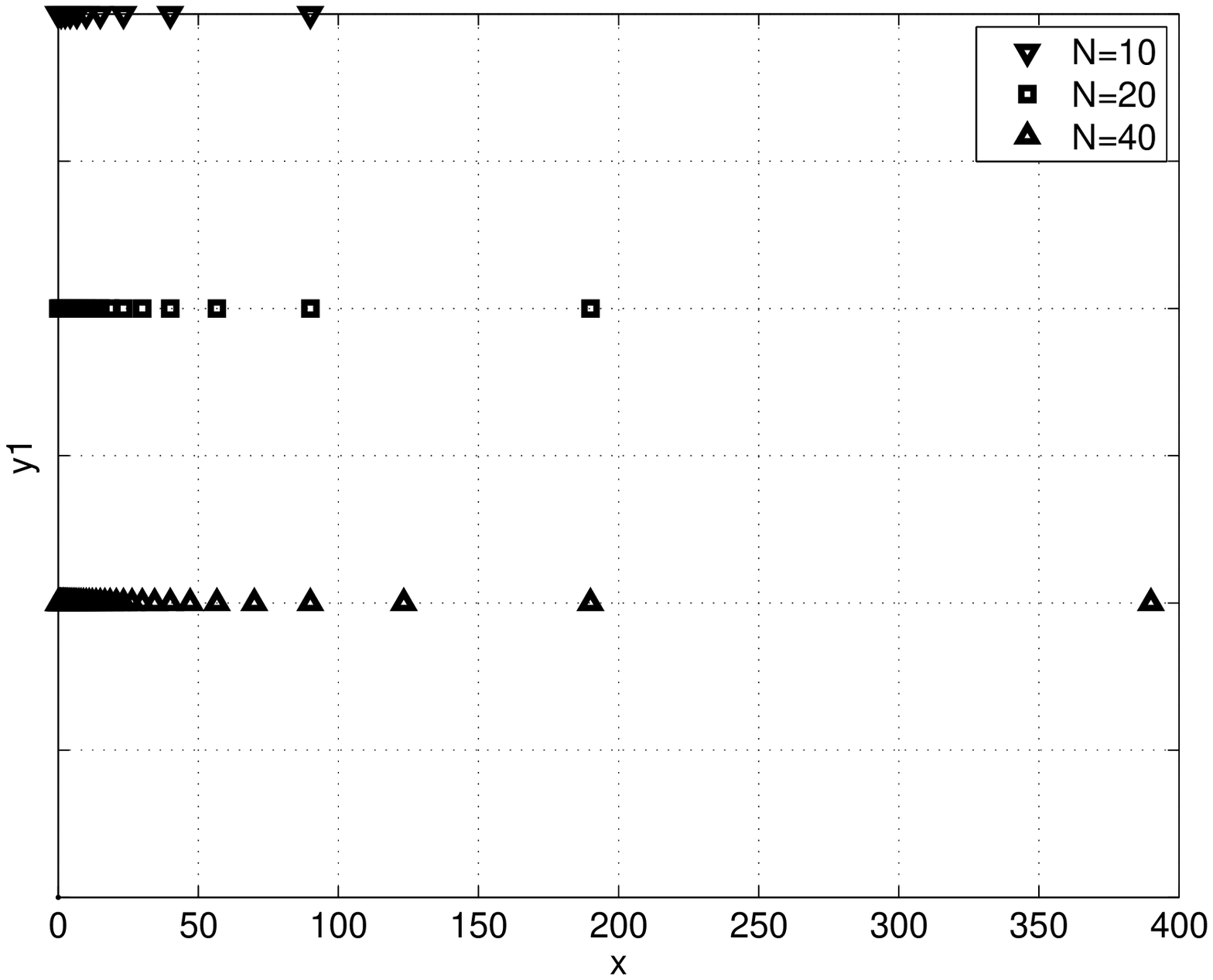}
\put(5,157){$x_N = \infty$} 
\put(5,108){$x_N = \infty$} 
\put(5,60){$x_N = \infty$}
\caption{\it Quasi-uniform grids: top frame for (\ref{eq:qu1}) and bottom frame for (\ref{eq:qu2}). We notice that, in both cases, the last mesh-point is $x_N = \infty$.}
\label{fig:m1N20}
\end{figure}

In order to derive the finite difference formulae, for the sake of simplicity, we consider a generic scalar variable $u(x)$.
We can approximate the values of this scalar variable at mid-points of the grid by
\begin{equation}
u_{n+1/2} \approx \frac{x_{n+3/4}-x_{n+1/2}}{x_{n+3/4}-x_{n+1/4}} u_n + \frac{x_{n+1/2}-x_{n+1/4}}{x_{n+3/4}-x_{n+1/4}} u_{n+1} \ .
\label{eq:u:mod}
\end{equation}
As far as the first derivative is concerned we can apply the following approximation
\begin{equation}
\left. \frac{du}{dx}\right|_{n+1/2} \approx \frac{u_{n+1}-u_n}{2\left(x_{n+3/4} - x_{n+1/4}\right)} \ .
\label{eq:du}
\end{equation}
These formulae use the value $ u_N = u(\infty) $, but not $ x_N = \infty $.
%The two finite difference approximations (\ref{eq:u:mod}) and (\ref{eq:du}) have order of accuracy $O(N^{-2})$. 
For a system of differential equations, (\ref{eq:u:mod}) and (\ref{eq:du}) can be applied component-wise.

\subsection{A non-standard finite difference scheme}\label{SS:NSFDS}
A non-standard finite difference scheme on a quasi-uniform grid for the class of BVPs
(\ref{p}) can be defined by using the approximations given by (\ref{eq:u:mod}) and (\ref{eq:du}) above.
A finite difference scheme for
(\ref{p}) can be written as follows:
\begin{eqnarray}
& {\bf U}_{n+1} - {\bf U}_{n} - a_{n+1/2} {\bf f} \left( x_{n+1/2}, b_{n+1/2}{\bf U}_{n+1} + c_{n+1/2}{\bf U}_{n} \right) = {\bf 0}
\ , \nonumber\\
& \mbox{for} \quad n=0, 1, \dots , N-1
\label{boxs} \\ 
& {\bf g} \left( {\bf U}_0,{\bf U}_N \right) = {\bf 0} \ ,  \nonumber
\end{eqnarray}
where 
\begin{eqnarray}\label{eq:abc}
a_{n+1/2} &=& 2\left(x_{n+3/4} - x_{n+1/4}\right) \ , \nonumber \\
b_{n+1/2} &=& \frac{x_{n+1/2}-x_{n+1/4}}{x_{n+3/4}-x_{n+1/4}} \ , \\
c_{n+1/2} &=&  \frac{x_{n+3/4}-x_{n+1/2}}{x_{n+3/4}-x_{n+1/4}} \nonumber \ , 
\end{eqnarray}
for $n=0, 1, \dots , N-1$.
The finite difference formulation (\ref{boxs}) has order of accuracy $O(N^{-2})$. It is evident that (\ref{boxs}) is a nonlinear system of $ d \; (N+1)$ equations in the $ d \; (N+1)$ unknowns $ {\bf U} = ({\bf U}_0,{\bf U}_1, \dots , {\bf U}_N)^T $. 
For the solution of (\ref{boxs}) we can apply the classical Newton's method along with the simple termination criterion
\begin{equation}\label{eq:Tcriterion}
{\ds \frac{1}{d(N+1)} \sum_{\ell =1}^{d} \sum_{n=0}^{N}
|\Delta {}^\ell U_{n}| \leq {\rm TOL}} \ ,
\end{equation}
where $ \Delta {}^ \ell U_{n} $, for $ n = 0,1, \dots, N $ and $ \ell = 1, 2, \dots , d $, is the difference between two successive iterate components and $ {\rm TOL} $ is a fixed tolerance. 

\section{\RE \ and error estimate}\label{S:extra}
The utilization of a quasi-uniform grid allows us to improve our numerical results.
The algorithm is based on Richardson's extrapolation, introduced by Richardson in \cite{Richardson:1910:DAL,Richardson:1927:DAL}, and it is the same for many finite difference methods: for numerical differentiation or integration, solving systems of ordinary or partial differential equations, see, for instance, \cite{Sidi:PEM:2003}.
To apply Richardson's extrapolation, we carry on several calculations on embedded uniform or quasi-uniform grids with total number of nodes $N_g$ for $g = 0$, $1$ , $\dots$, $G$: e.g., for the numerical results reported in the next section we have used $5$, $10$, $20$, $40$, $80$, $160$, $320$, $640$, $1280$, $2560$, or $5120$ grid-points. 
We can identify these grids with the index $g=0$, the coarsest one, $1$, $2$, and so on towards the finest grid denoted by $g = G$.
Between two adjacent grids, all nodes of largest steps are identical to even nodes of denser grid due to the uniformity. 
To find a more accurate approximation we can apply $k$ Richardson's extrapolations on the used grids
\begin{equation}\label{eq:Rextra}
U_{g+1,k+1} = U_{g+1,k} + \frac{U_{g+1,k}-U_{g,k}}{2^{p_k}-1} \ ,
\end{equation}
where $g \in \{0, 1, 2 , \dots , G-1\}$, $k \in \{0, 1, 2, \dots , G-1\}$, $2 = N_{g+1}/N_{g}$ appearing in the denominator is the grid refinement ratio, and $p_k$ is the true order of the discretization error.
%This formula is asymptotically exact in the limit as $N$ goes to infinity if we use uniform or quasi-uniform grids.
We notice that to obtain each value of $U_{g+1,k+1}$ requires having computed two solution $U$ in two embedded grids, namely $g+1$ and $g$ at the extrapolation level $k$.
For any $g$, the level $k=0$ represents the numerical solution of $U$ without any extrapolation, which is obtained as described in subsection \ref{SS:NSFDS}.
In this case, Richardson extrapolation uses two solutions on embedded refined grids to define a more accurate solution that is reliable only when the grids are sufficiently fine. 
The case $k=1$ is the classical single Richardson's extrapolation, which is usually used to estimate the discretization error or to improve the solution accuracy.
If we have computed the numerical solution on $G+1$ nested grids then we can apply equation (\ref{eq:Rextra}) $G$ times performing $G$ Richardson's extrapolations. 

The theoretical orders $p_k$ of accuracy of the numerical solution $U$ with $k$ extrapolations verify the relation
\begin{equation}\label{eq:pk}
p_k = p_0 + k (p_1-p_0) \ ,
\end{equation}
where this equation is valid for $k \in \{0, 1, 2, \dots , G-1\}$.
In any case, the values of $p_k$ can be obtained a priori by using appropriate Taylor series or a posteriori by
\begin{equation}\label{eq:pk:calc}
p_k \approx {\ds \frac{\log(|U_{g,k}-u|)-\log(|U_{g+1,k}-u|)}{\log(2)}} \ ,
\end{equation}
where $u$ is again the exact solution (or, if the exact solution is unknown, a reference solution computed with a suitable large value of $N$) evaluated at the same grid-points of the numerical solution. 

%\subsection{Error estimate}\label{SS:error}
To show how Richardson's extrapolation can be also used to get an error estimate for the computed numerical solution we use two numerical solutions $U_{N}$ and $U_{2N}$ computed by doubling the number of grid-points.
Taking into account equation (\ref{eq:Rextra}) we can conclude that the error estimate by Richardson's extrapolation is given by
\begin{equation}\label{eq:est1}
E = \frac{U_{2N}-U_{N}}{2^{p_0}-1} \ ,
\end{equation}
where $p_0$ is the true order of the discretization error.
Hence, $E$ is an estimation of the truncation error found without knowledge of the exact solution.
We notice that $E$ is the error estimate for the more accurate numerical solution $U_{2N}$ but only on the grid points of $U_N$. 
%Moreover, using the above mentioned grids we have to replace $r = 2$ in equations (\ref{eq:est1}) and (\ref{eq:p}).

\section{Numerical results: a BVP in colloids theory}
In this section, we consider a benchmark problem with known exact solution for our error estimator.
It should be mentioned that all numerical results reported in this paper were performed on an ASUS personal computer with i7 quad-core Intel processor and 16 GB of RAM memory running Windows 8.1 operating system. 
The non-standard finite difference scheme described above has been implemented in FORTRAN.
The numerical results reported in this section were computed by setting 
\begin{equation}\label{eq:TOL}
{\rm TOL} = 1\E-12 \ .
\end{equation}

The benchmark problem, see Alexander and Johnson \cite{Alexander:1949:CS}, arises within the theory of colloids and is given by 
\begin{align}\label{colloid:model}
& {\ds \frac{d^2 u}{dx^2}} - 2 \sinh(u) = 0 \qquad \ x \in [0, \infty] \ , \nonumber \\[-1ex]
& \\[-1ex]
& u(0) = u_0 \ , \qquad u(\infty) = 0 \ , \nonumber 
\end{align}
where $u_0 > 0$.
The exact solution of the BVP (\ref{colloid:model}) 
\begin{equation}\label{eq:colloid:model:exact}
u(x) = 2 \; \ln\left(\frac{(e^{u_0/2}+1)\; e^{\sqrt{2}\; x}+(e^{u_0/2}-1)}{(e^{u_0/2}+1)\; e^{\sqrt{2}\; x}-(e^{u_0/2}-1)}\right) \ ,
\end{equation}
has been found by Countryman and Kannan \cite{Countryman:1994:CNB,Countryman:1994:NBV}, and the missing initial condition is given by
\begin{equation}\label{eq:colloid:model:exact:dudx0}
\frac{du}{dx}(0) = -2\; \sqrt{\cosh(u_0) -1} \ .
\end{equation}

\begin{figure}[!hbt]
\centering
\psfrag{x}[][]{\small $x$} 
\psfrag{U}[][]{\tiny ${}^1U$} 
\psfrag{dU}[][]{\tiny ${}^2U$} 
\psfrag{u}[][]{\tiny ${}^1u$} 
\psfrag{du}[][]{\tiny ${}^2u$} 
\psfrag{UdU}[][]{\small ${}^1U$, ${}^2U$} 
\psfrag{udu}[][]{\small ${}^1U$, ${}^2U$, ${}^1u$, ${}^2u$} 
\includegraphics[width=\www]{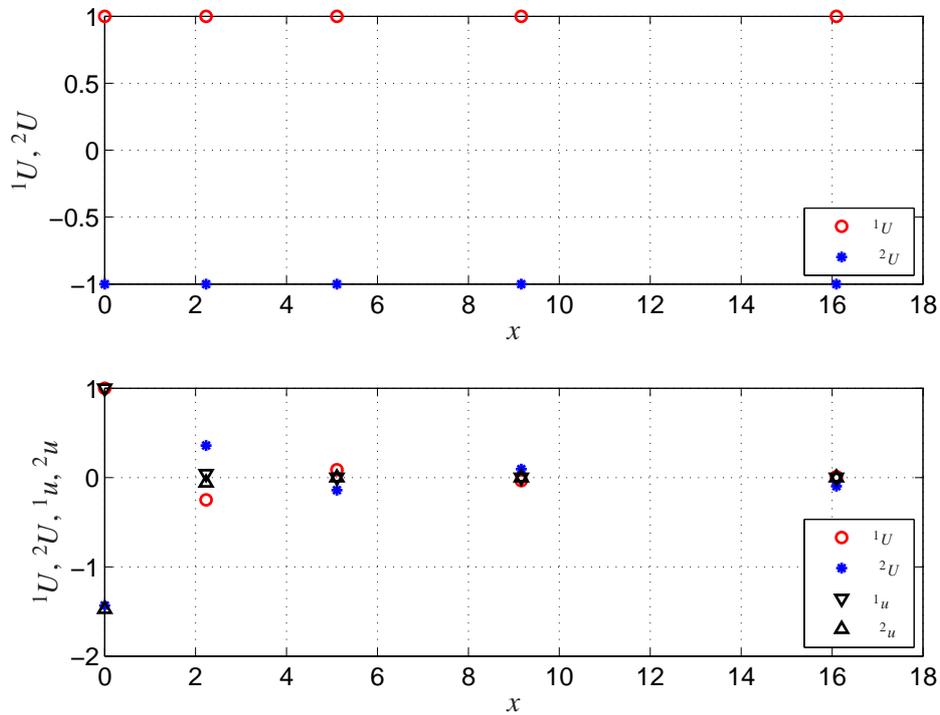}
\caption{\it Sample iterates for problem (\ref{colloid:model}) with $u_0 = 1$.}
\label{fig:it}
\end{figure}

\begin{figure}[!hbt]
\centering
\psfrag{x}[][]{\small $x$} 
\psfrag{udu}[][]{\small ${}^1U$, ${}^2U$, ${}^1u$, ${}^2u$} 
\psfrag{U}[][]{\tiny ${}^1U$} 
\psfrag{dU}[][]{\tiny ${}^2U$} 
\psfrag{u}[][]{\tiny ${}^1u$} 
\psfrag{du}[][]{\tiny ${}^2u$} 
\includegraphics[width=\www]{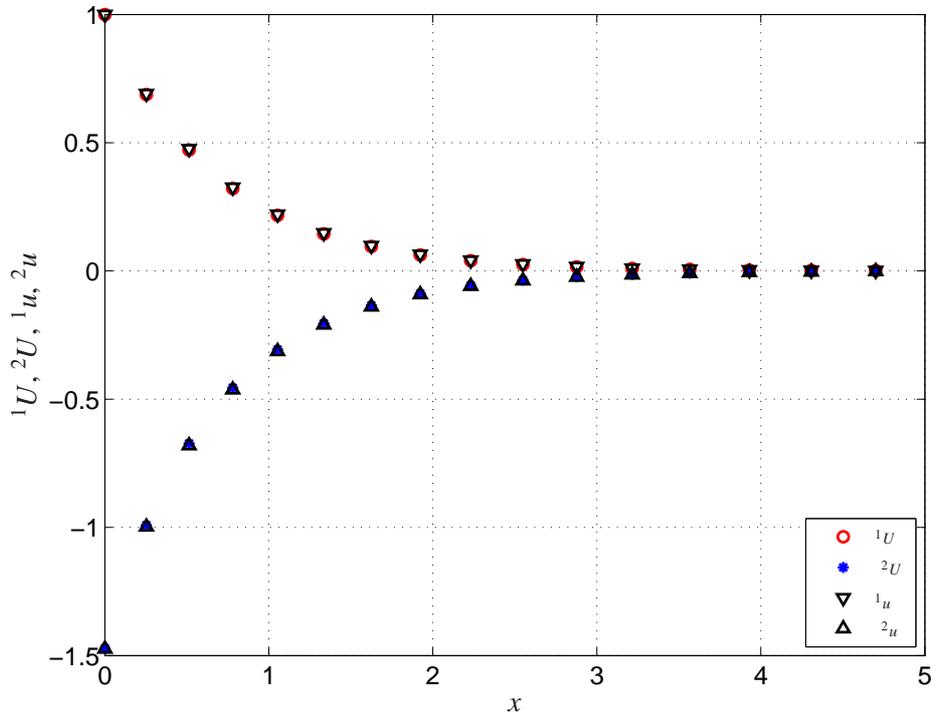}
\caption{\it Final iterate and exact solution for problem (\ref{colloid:model}) with $u_0 = 1$ and $N=40$: zoom of the transitory region. %On the left: solutions on the quasi-uniform grid and on the right: zoom of the transitory region.
}
\label{fig:sol_1}
\end{figure}

We rewrite the governing differential equation as a first order system and indicate the exact solution with ${\bf u} = ({}^1u, {}^2u)^T$ and the numerical solution with ${\bf U} = ({}^1U, {}^2U)^T$.
In order to fix a specific problem, as a first test case, we consider $u_0=1$.
As mentioned before we used $5$, $10$, $20$, $40$, $80$, $160$, $320$, $640$, $1280$, $2560$, or $5120$ grid-points, so that $G=10$, and we adopted a continuation approach for the choice of the first iterate.
This means that the accepted solution for $N=5$ is used as first iterate for $N=10$, where the new grid values are approximated by linear interpolations, and so on.
The first iterate for the grid with $N=5$, where the field variable was taken constant and equal to one and its derivative was taken also constant and equal to minus one, is shown in the top frame of figure \ref{fig:it}.
The bottom frame of the same figure shows the accepted numerical solution.
Our relaxation algorithm takes seven iterations to verify the termination criterion (\ref{eq:Tcriterion}) with TOL given by (\ref{eq:TOL}).
Once the continuation approach has been initialized, the iteration routine needs 3 or 4 iterations to get a numerical solution that verifies the stopping criterion. 
For the sake of completeness in figure \ref{fig:sol_1} we display the numerical solution for $N=40$ along with the exact solution. 

\begin{figure}[!hbt]
\centering
\psfrag{ERR}[][]{$e$}% \small $E_{g,k}$ for $k=0,1,2$} 
\psfrag{N}[l][]{\small $ \quad N_g$} 
\includegraphics[width=.45\textwidth]{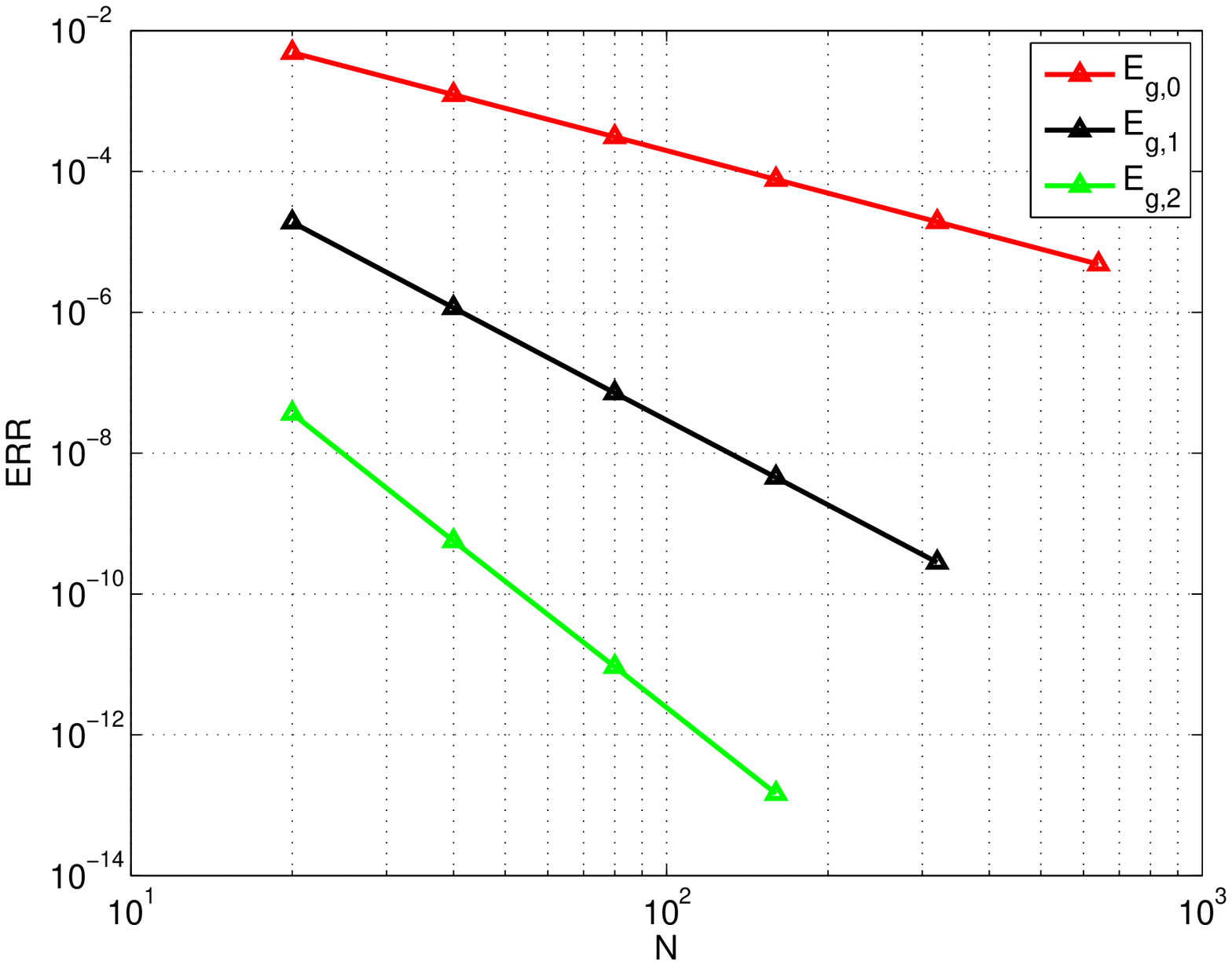} \hfill
\includegraphics[width=.45\textwidth]{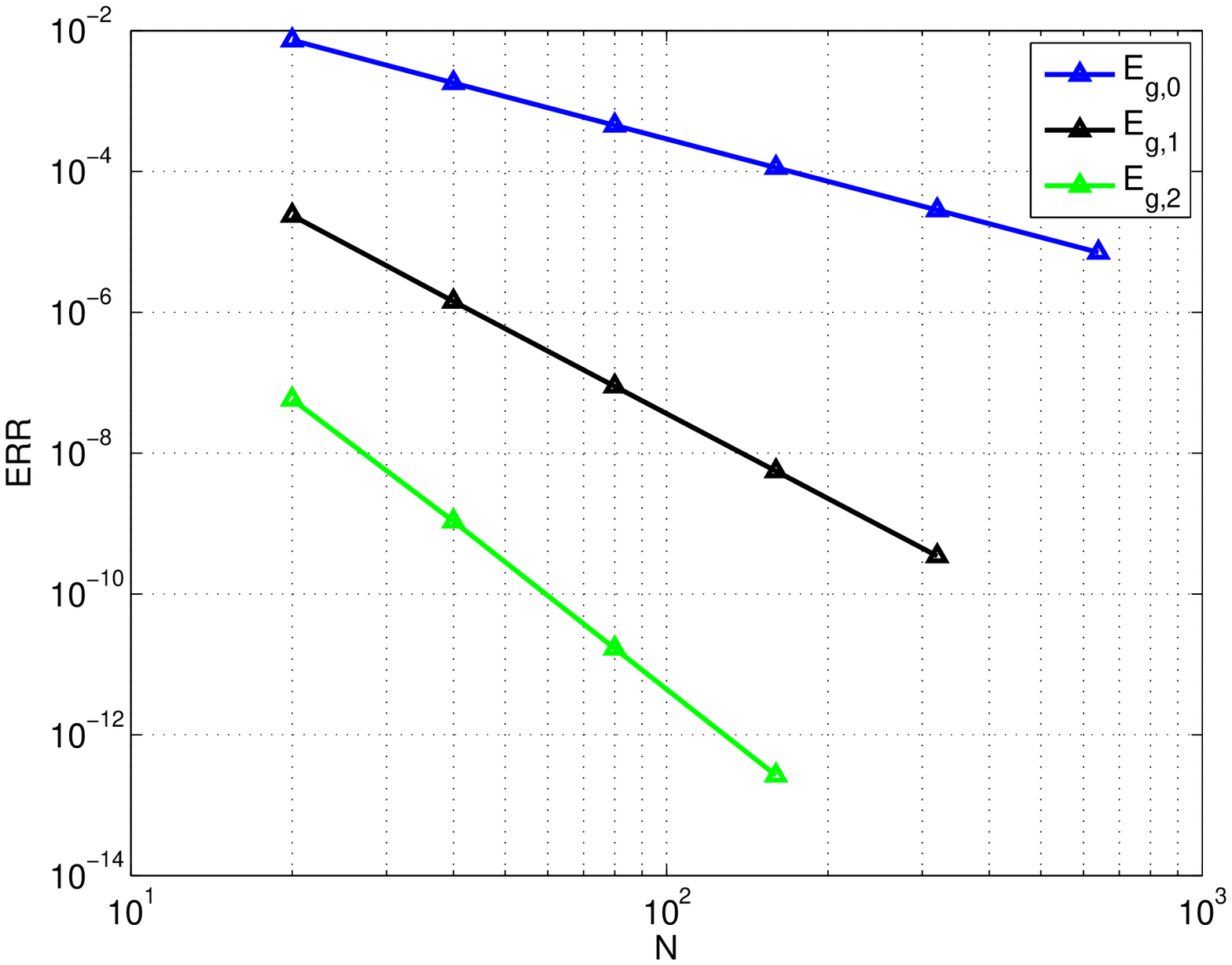}
\caption{\it Graphical derivation of the values of $p_0$, $p_1$ and $p_2$ for the non standard finite difference scheme applied to (\ref{colloid:model}) with $u_0 = 1$.}
\label{fig:Eforp}
\end{figure}

\begin{figure}[!hbt]
\centering
\psfrag{x}[][]{\small $x$} 
\psfrag{e}[][]{\small ${}^1E$, ${}^2E$, ${}^1e$, ${}^2e$} 
\psfrag{EU}[][]{\tiny ${}^1E$} 
\psfrag{EdU}[][]{\tiny ${}^2E$} 
\psfrag{eU}[][]{\tiny ${}^1e$} 
\psfrag{edU}[][]{\tiny ${}^2e$} 
\includegraphics[width=.7\textwidth]{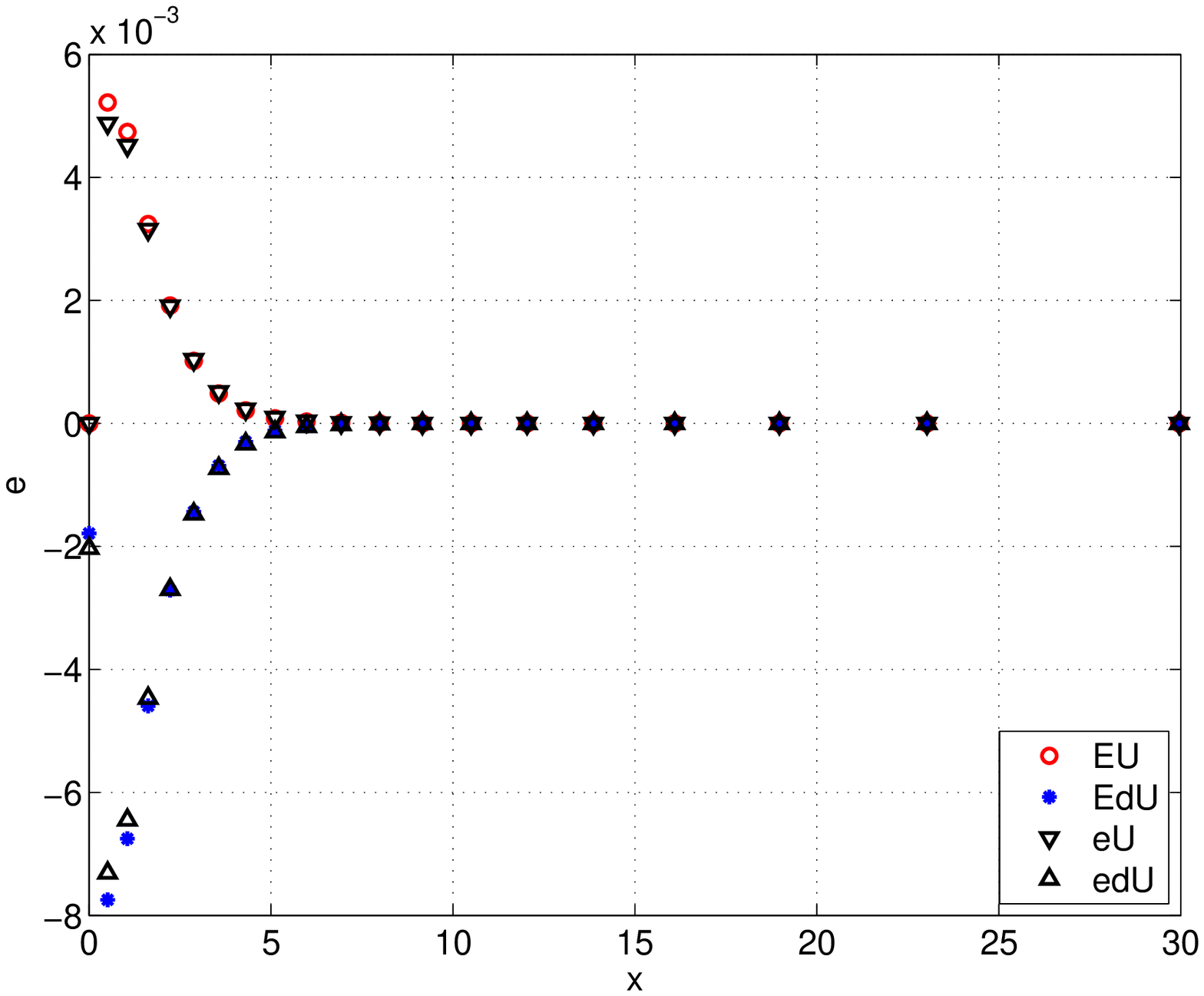} \\
\includegraphics[width=.7\textwidth]{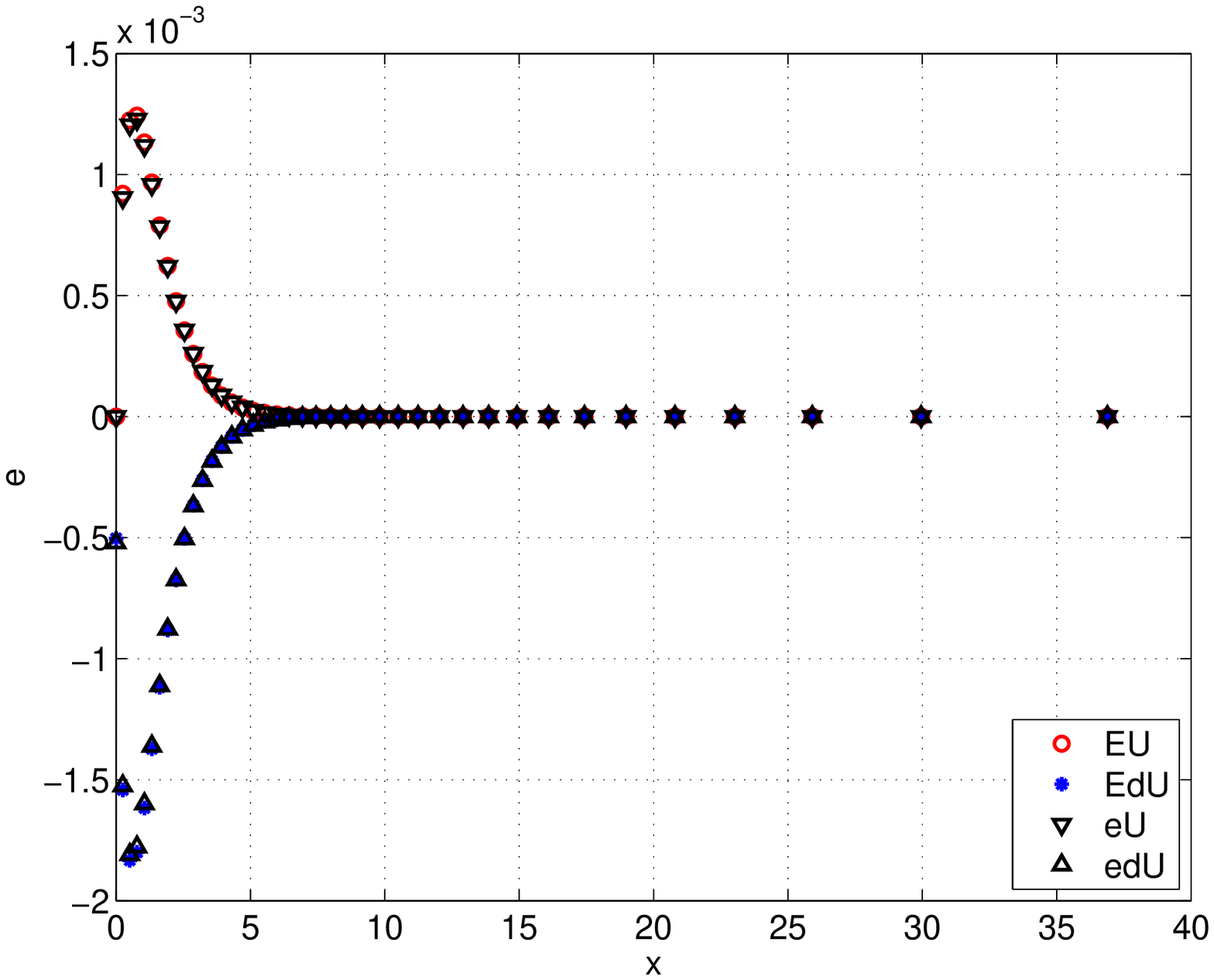}
\caption{\it Global and a posteriori error estimates for the field variable and its first derivative for (\ref{colloid:model}) with $u_0 = 1$. Top: $N = 20$, $40$, and bottom: $N = 40$, $80$. Here ${}^1e$ and ${}^2e$ are the global errors by equation (\ref{eq:GE}), whereas ${}^1E$ and ${}^2E$ are the error estimates provided by equation (\ref{eq:est1}).}
\label{fig:Errors_1}
\end{figure}

Figure \ref{fig:Eforp} shows in a log by log scale the computed errors.
We can compute the order of accuracy $p_0$, $p_1$ and $p_2$ according to the formula (\ref{eq:pk:calc}).
As it is easily seen from figure \ref{fig:Eforp} we got $p_0 \approx 2$, $p_1 \approx 4$ and $p_2 \approx 6$ for both the field variable and its first derivative.
As far as the a posteriori error estimator is concerned in figure \ref{fig:Errors_1} we report the computation related to two sample cases: namely, the estimate obtained by using $N = 20$, $40$ and $N = 40$, $80$. 
We notice that the global error, for both the solution components, is of order $10^{-3}$ and it decreases as we refine the grid. 
It is easily seen that the estimator defined by equation (\ref{eq:est1}) provides upper bounds for the global error.

\begin{figure}[!hbt]
\centering
\psfrag{x}[][]{\small $x$} 
\psfrag{udu}[][]{\small ${}^1U$, ${}^2U$, ${}^1u$, ${}^2u$} 
\psfrag{U}[][]{\tiny ${}^1U$} 
\psfrag{dU}[][]{\tiny ${}^2U$} 
\psfrag{u}[][]{\tiny ${}^1u$} 
\psfrag{du}[][]{\tiny ${}^2u$} 
\includegraphics[width=\www]{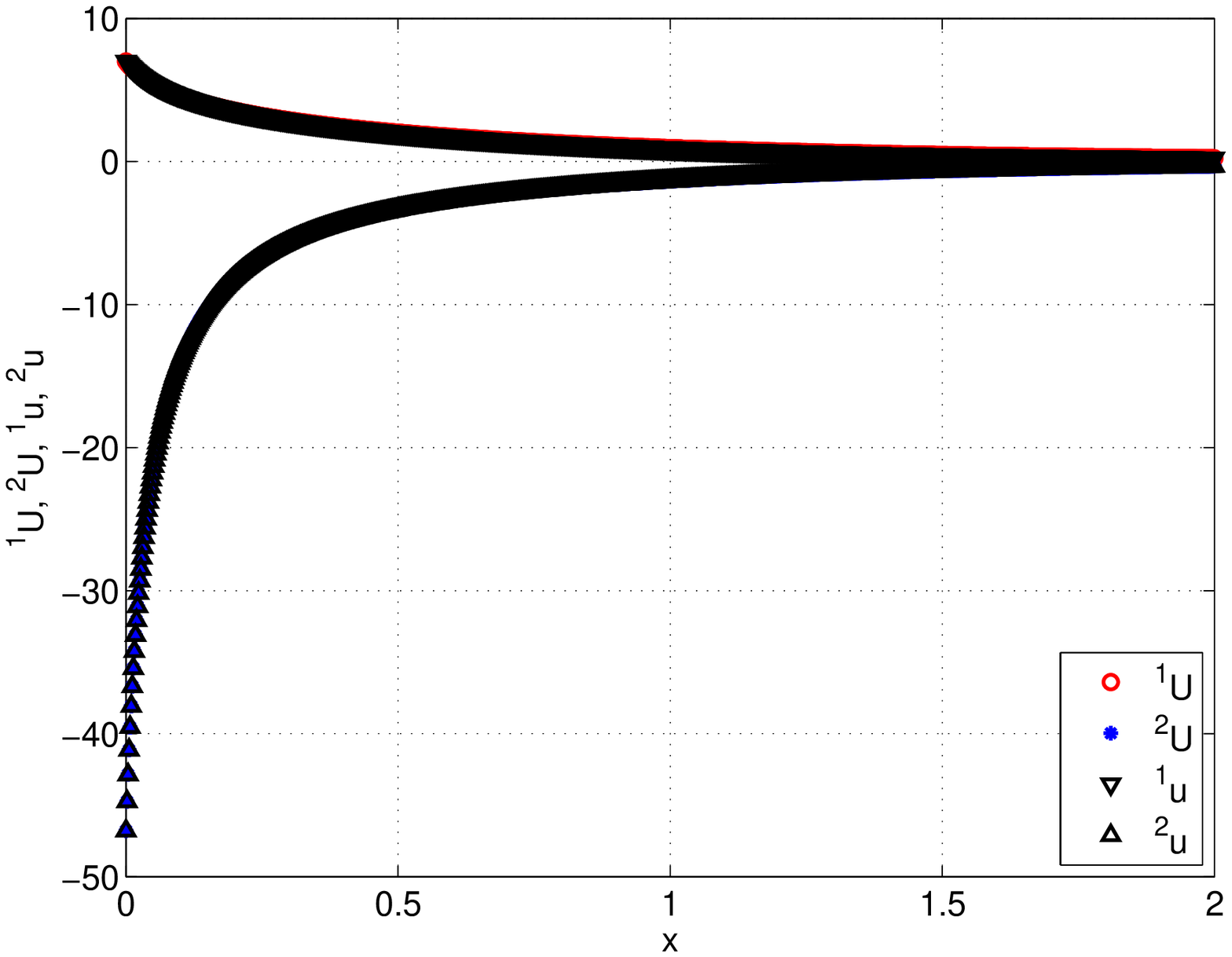}
\caption{\it Final iterate and exact solution for problem (\ref{colloid:model}) with $u_0 = 7$ and $N=5120$: zoom of the transitory region. %On the top: solutions on the quasi-uniform grid and on the bottom: zoom of the transitory region.
}
\label{fig:sol_7}
\end{figure}

A more challenging test case is given by setting $u_0=7$.
In figure \ref{fig:sol_7} we display the numerical solution for $N=5120$ along with the exact solution. 
%\textcolor{red}{NB: Tabella per il valore di $du/dx(0)$.}
\begin{table}[!hbt]
\renewcommand\arraystretch{1.4}
	\centering
\begin{tabular}{rccc}
\hline \\[-1.5ex]
$N_g$ & ${}^2U_{g,0}$ & ${}^2U_{g,1}$ & ${}^2U_{g,2}$ \\[1.5ex]
\hline
%20    & $-24.037154011637909$ & $ $ & $ $ \\
%40    & $-32.118139297376729$ &  $-34.811801059289671$ & $ $\\
%80    & $-39.285772793638444$ & $-41.674983959059013$ & $ $ \\ 
160   & $-43.835177171609345$ & $$ & $$ \\  
320   & $-45.864298511341850$ & $-46.540672291252690$ & $$ \\  
640   & $-46.537797149336093$ & $-46.762296695334179$ & $-46.777071655606278$ \\  
1280   & $-46.725033491934731$ & $-46.787445606134277$ & $-46.789122200187620$ \\  
2560   & $-46.773360098843838$ & $-46.789468967813541$ & $-46.789603858592159$ \\  
5120   & $-46.785544794016836$ & $-46.789606359074504$ & $-46.789615518491907$ \\  
\hline
\end {tabular}
\caption{\it Richardson's extrapolation for $\frac{du}{dx}(0)={}^2U_0$.}
\label{tab:extra} 
\end{table}
In table \ref{tab:extra} we list the computed as well as the extrapolated values obtained for the missing initial condition. 
For the sake of brevity, in this table, we do not report the fewer accurate values obtained with the coarser grids.
These results can be compared with the exact value, $\frac{du}{dx}(0) \approx -46.789615734913319$, obtained by equation (\ref{eq:colloid:model:exact:dudx0}).

\begin{figure}[!hbt]
\centering
\psfrag{E0}[][]{\small $E_{g,0}$} 
\psfrag{ERR}[][]{$e$}%\small $E_{g,k}$ for $k=0,1,2$} 
\psfrag{N}[][]{\small $N_g$} 
\includegraphics[width=.45\textwidth]{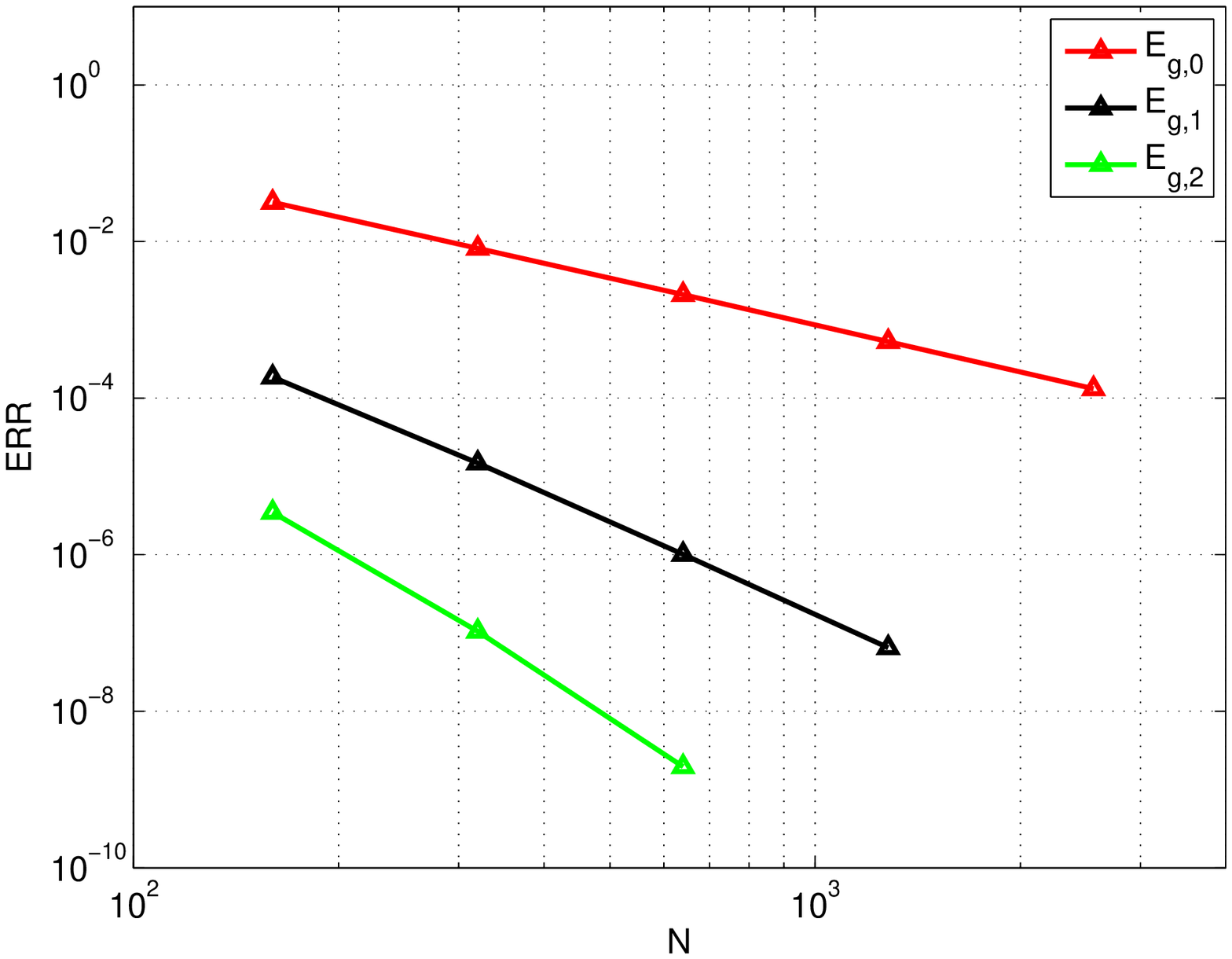} \hfill
\includegraphics[width=.45\textwidth]{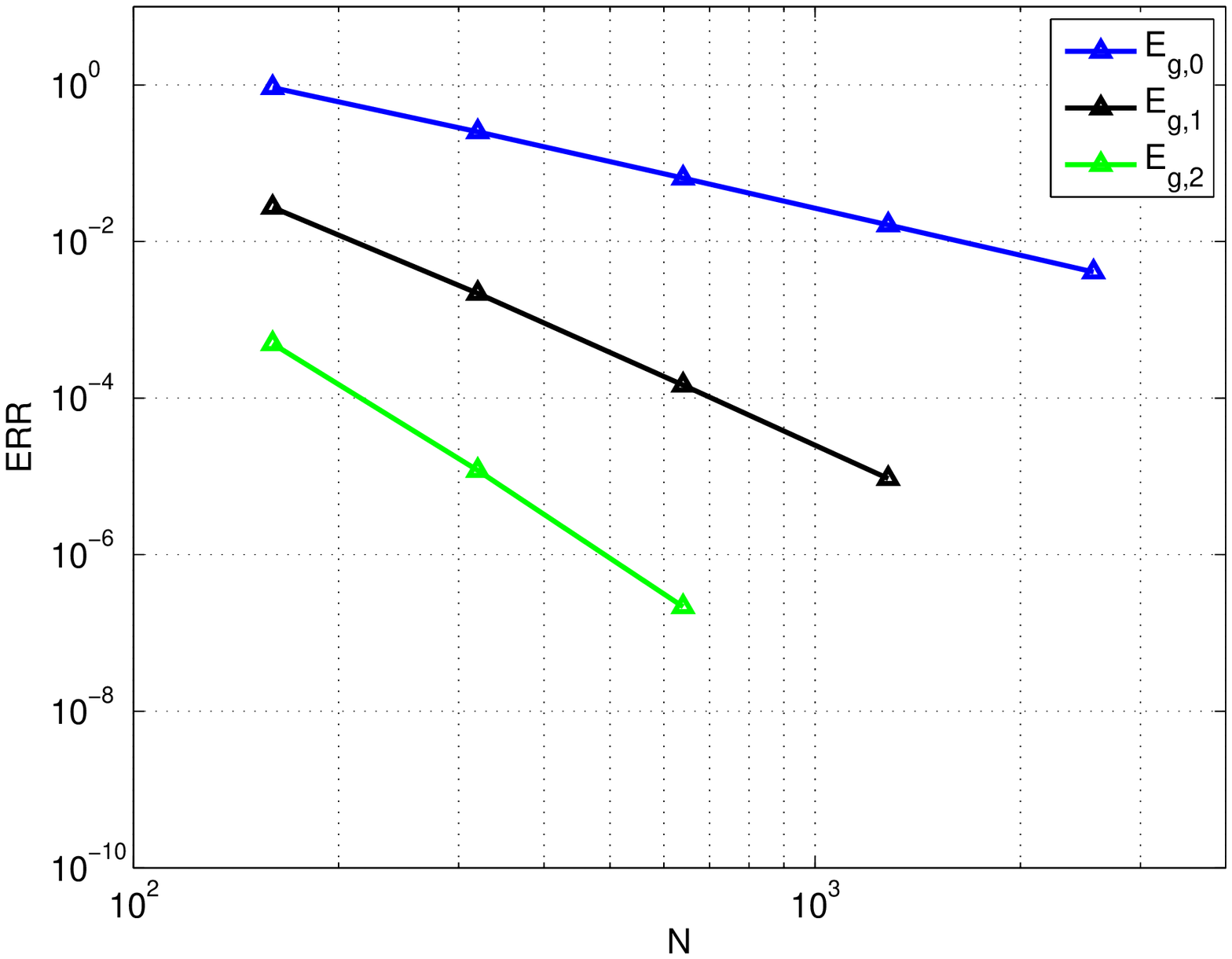}
\caption{\it Problem (\ref{colloid:model}) with $u_0 = 7$, global errors for: the field variable on the left and its first derivative on the right.}
\label{fig:Errors:c7}
\end{figure}

Figure \ref{fig:Errors:c7} shows in a log by log scale the computed errors.
From figure \ref{fig:Errors:c7} it is clear that the computed orders, using equation (\ref{eq:pk:calc}), are slightly different from the theoretical ones, namely: $p_0 \approx 1.99$, $p_1 \approx 3.96$ and $p_2 \approx 5.77$.

\begin{figure}[!hbt]
\centering
\psfrag{x}[][]{\small $x$} 
%\psfrag{e}[][]{\small ${}^1E$, ${}^2E$, ${}^1e$, ${}^2e$} 
%\psfrag{EU}[][]{\tiny ${}^1E$} 
%\psfrag{EdU}[][]{\tiny ${}^2E$} 
%\psfrag{eU}[][]{\tiny ${}^1e$} 
%\psfrag{edU}[][]{\tiny ${}^2e$} 
\includegraphics[width=.7\textwidth]{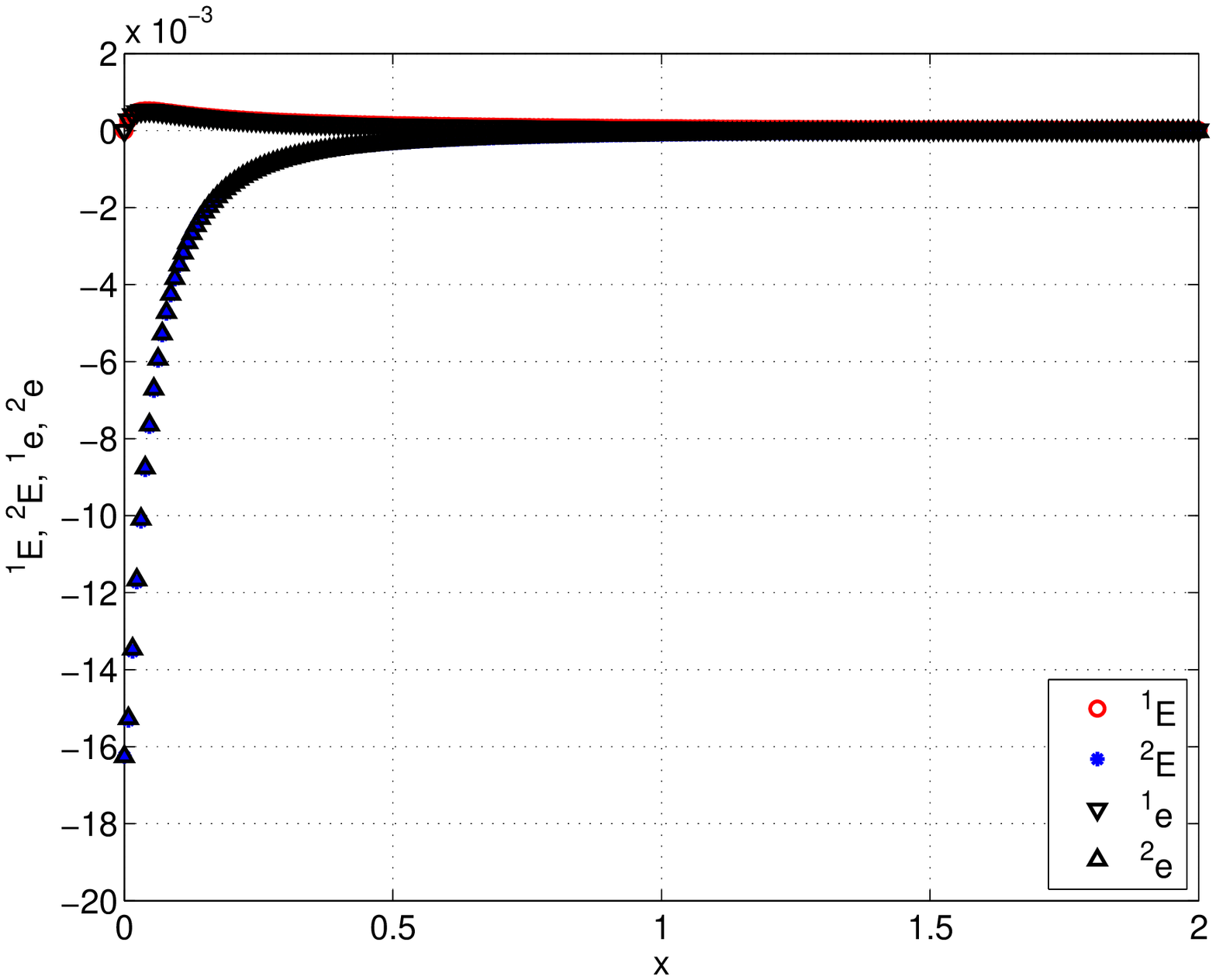} \\
\includegraphics[width=.7\textwidth]{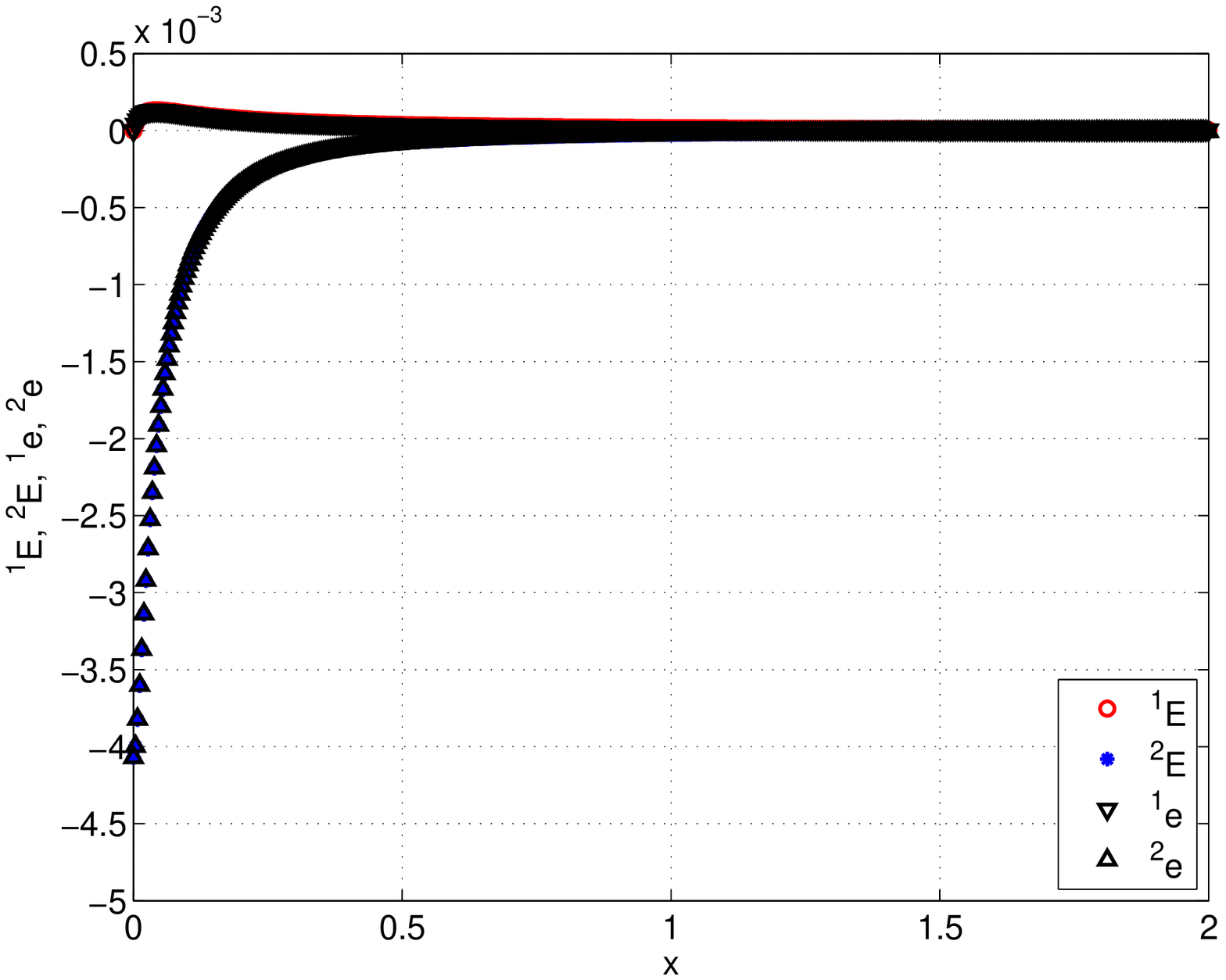}
\caption{\it Zoom in the domain related to the initial transitory for global errors and a posteriori error estimates for the field variable and its first derivative for (\ref{colloid:model}) with $u_0 = 7$. Top:  $N = 1280$, $2560$ and bottom: $N = 2560$, $5120$. Here ${}^1e$ and ${}^2e$ are the global errors by equation (\ref{eq:GE}), whereas ${}^1E$ and ${}^2E$ are the error estimates provided by equation (\ref{eq:est1}).}
\label{fig:Errors_7}
\end{figure}
As far as the a posteriori error estimator is concerned in figure \ref{fig:Errors_7} we report the computation related to two sample cases: namely, the estimate obtained by using $N = 1280$, $2560$ and $N = 2560$, $5120$. 
Once again the global error, for both the solution components, decreases as we refine the grid and the estimator defined by equation (\ref{eq:est1}) provides upper bounds for the global error.

%In the following 
\section{Conclusions}
In this paper, we have defined a posteriori estimator for the global error of a non-standard finite difference scheme applied to boundary value problems defined on an infinite interval.
A test problem was examined for which the exact solution is known and we tested our error estimator for two sample cases: a simpler one with smooth solutions component and a more challenging one presenting an initial fast transitory for one of the solution components.
For this second test case, we showed how Richardson extrapolation can be used to improve the numerical solution using the order of accuracy and numerical solutions from two nested quasi-uniform grids. 
Moreover, the reported numerical results clearly show that our non-standard finite difference scheme, implemented along with the error estimator defined in this work, can be used to solve challenging problems arising in the applied sciences. 

In our previous paper \cite{Fazio:2014:FDS} we derived instead of equation (\ref{eq:u:mod}) the finite difference formula
\begin{equation}
u_{n+1/2} \approx \frac{x_{n+1}-x_{n+1/2}}{x_{n+1}-x_n} u_n + \frac{x_{n+1/2}-x_{n}}{x_{n+1}-x_n} u_{n+1} \ .
\label{eq:u}
\end{equation}
However, by setting $n =N-1$ this formula (\ref{eq:u}), replacing $x_N = \infty$, reduces to $u_{N-1/2} = u_{N-1}$ that does not involve the boundary value $u_N$ and therefore the boundary condition cannot be used.
In that work, this was the reason that forced us to modify this formula at $ n = N-1$, see \cite{Fazio:2014:FDS} for details.
The mentioned drawback is completely overcome by the new formula (\ref{eq:u:mod}).

\vspace{1.3cm}

\noindent {\bf Acknowledgement.} {The research of this work was 
supported, in part, by the University of Messina and by the GNCS of INDAM.}
 
\bibliographystyle{plain}
\bibliography{/ricerca/latex2e/bib/fazio200,/ricerca/latex2e/bib/blt,/ricerca/latex2e/bib/nppud,/ricerca/latex2e/bib/bvpud,/ricerca/latex2e/bib/sc-gita,/ricerca/latex2e/bib/sc-gita2,/ricerca/latex2e/bib/na,/ricerca/latex2e/bib/nahistory,/ricerca/latex2e/bib/Math,/ricerca/latex2e/bib/FreeBVP,/ricerca/latex2e/bib/mbpp}

\end{document}